\documentclass[11pt,a4paper]{amsart}
\usepackage{color}
\usepackage{amsmath,amssymb,latexsym,amsthm}
\usepackage[utf8]{inputenc}

\usepackage{tikz}
\usetikzlibrary{arrows}

\newtheorem{Thm}{Theorem}[section]
\newtheorem{Lemma}[Thm]{Lemma}
\newtheorem{Prop}[Thm]{Proposition}

\newtheorem{Cor}[Thm]{Corollary}
\newtheorem{Rem}[Thm]{Remark}

\newtheorem{Example}[Thm]{Example}

\newtheorem{Question}{Question}

\newcommand{\bb}{\beta}

\newcommand{\e}{\varepsilon}

\newcommand{\la}{\lambda}

\newcommand{\s}{\kappa}
\newcommand{\G}{\Gamma}
\newcommand{\D}{\Delta}

\newcommand{\dist}{\textrm{dist}}

\newcommand{\Lc}{\textrm{L}}
\newcommand{\dns}{\textrm{DN-S}}
\newcommand{\dss}{\textrm{DSS}}

\newcommand{\ds}{\textrm{d}}

\newcommand{\N}{\mathbb{N}}
\newcommand{\R}{\mathbb{R}}

\begin{document}

\title[A quantitative approach to $\dns$ operators]{A quantitative approach to disjointly\\ non-singular operators}

\thanks{Supported in part by MICINN (Spain), Grant PID2019-103961GB-C22.
}

\author[M.\ Gonz\'alez]{Manuel Gonz\'alez}
\address{Departamento de Matem\'aticas, Facultad de Ciencias,
Universidad de Canta\-bria, E-39071 Santander, Spain}
\email{manuel.gonzalez@unican.es}

\author[A.\ Martin\'on]{Antonio Martin\'on}
\address{Departamento de An\'alisis Matem\'atico, Facultad de Ciencias,
Universidad de La Laguna, E-38271 La Laguna (Tenerife), Spain} \email{anmarce@ull.es }


\begin{abstract}
We introduce and study some operational quantities which characterize the disjointly non-singular operators from a Banach lattice $E$ to a Banach space $Y$ when $E$ is order continuous, and some other quantities which characterize the disjointly strictly singular operators for arbitrary $E$. 
\end{abstract}

\maketitle

\thispagestyle{empty}

\section{Introduction}

The disjointly strictly singular operators ($\dss$ operators) were introduced in \cite{Hernandez-Salinas:89} as those operators $T:E\to Y$ from a Banach lattice $E$ into a Banach space $Y$ such that $T$ is not an isomorphism in any  subspace of $E$ generated by a disjoint sequence of non-zero vectors. These operators have been useful in the study of the structure of Banach lattices (see \cite{FHT:16}, \cite{FLT:16} and references therein). 
More recently, the disjointly non-singular operators ($\dns$ operators) where introduced in \cite{GMM:20} (see also \cite{Bilokopytov:21}) as those operators  $T:E\to Y$ that are not strictly singular in any subspace of $E$ generated by a disjoint sequence of non-zero vectors. Note that the properties in the definition of these two classes are opposite.  

In this paper we study the classes of operators $\dss$ and $\dns$ from a quantitative point of view by introducing four operational quantities $\G_d(T)$, $\D_d(T)$, $\tau_d(T)$ and $\s_d(T)$. When $E$ is order continuous, $T\in\dns(E,Y)$ is equivalent to  $\G_d(T)>0$, or $\s_d(T)>0$; and for $E$ arbitrary, $T\in\dss(E,Y)$ is equivalent to $\D_d(T)=0$, or $\tau_d(T)=0$. 
These four quantities are inspired by some others introduced by Schechter \cite{schechter} in his study of Fredholm theory. 

In \cite{GMM:20}, the quantity $\beta(T) = \inf_{(x_n)} \liminf_{n\to\infty} \|Tx_n\|$, where the infimum is taken over the normalized disjoint sequences $(x_n)$ in $E$, was defined. We show that $T\in\dns(E,Y)$ if and only if $\beta(T)>0$ when $E$ is order continuous. This result was proved in \cite[Theorem 5.7]{Bilokopytov:21} using different techniques. We also prove that $\beta(T)\leq \G_d(T)$, but there is no $C>0$ such that $\G_d(T) \leq C\beta(T)$ for each $T\in\Lc(\ell_2,Y)$; hence $\G_d$ and $\beta$ are not equivalent. Moreover, $\tau_d(T)\leq \D_d(T)$, but the quantities $\tau_d$ and $\D_d$ are not equivalent.  

We also prove some inequalities for  these  operational quantities; e.g., for $T,S\in \Lc(E,Y)$, we have $\G_d(T+S)\leq \G_d(T)+\D_d(S)$. When $E$ is order continuous, this inequality allows us to improve the stability result for $\dns$ operators under $\dss$ perturbations  obtained in \cite{GMM:20}. 

\subsection*{Notation}
Throughout the paper $X$ and $Y$ are Banach spaces, and $E$ is a Banach lattice. The unit sphere of $X$ is  $S_X=\{x\in X : \|x\|=1\}$, and for a sequence $(x_n)$ in $X$, $[x_n]$ denotes the closed subspace generated by $(x_n)$.

All the operators are linear and bounded, and $\Lc(X,Y)$ denotes the set of all the 
operators from $X$ into $Y$. Given $T\in \Lc(X,Y)$, its {\it injection modulus} is  $j(T) := \inf_{\| x\| =1}  \|Tx \|$. Recall that $j(T)>0$ if and only if $T$ is an isomorphism from $X$ onto $TX$. We denote by $T_M$ the restriction of $T\in \Lc(X,Y)$ to a closed subspace $M$ of $X$. 

If $(\Omega,\Sigma,\mu)$ is a measure space, the domain of a measurable function $f: \Omega\to \R$ is the set $D(f) = \{ t\in \Omega : f(t) \neq 0 \}$, and $1_A$ denotes the characteristic function of $A\in\Sigma$. We write $L_p$ for $L_p[0,1]$, $1 \leq  p \leq \infty$.

\section{Preliminaries}

An operator $T\in\Lc(X,Y)$ is {\it strictly singular} if there is no closed infinite dimensional subspace $M$ of $X$ such that the restriction $T_M$ is an isomorphism, and $T$ is {\it upper semi-Fredholm} if its kernel is finite dimensional and its range is closed. 

An operator $T\in\Lc(E,Y)$ is {\it disjointly strictly singular} if there is no disjoint sequence of non-zero vectors $(x_n)$ in $E$ such that $T_{[x_n]}$ is an isomorphism. We denote by $\dss (E,Y)$ the set of all $T \in \Lc(E,Y)$ which are disjointly strictly singular. The class $\dss$ was introduced by Hernández and Rodríguez-Salinas in \cite{Hernandez-Salinas:89}. More information on this class can be found in \cite{Hernandez:90}.
\medskip

An operator $T\in \Lc(E,Y)$ is {\it disjointly non-singular} if there is no disjoint sequence of non-zero vectors $(x_n)$ in $E$ such that $T_{[x_n]}$ is strictly singular.  
We denote $\dns(E,Y)$ the set of all  $T \in \Lc(E,Y)$ which are disjointly non-singular. These operators were recently introduced in \cite{GMM:20}, and have been studied by Bilokopytov in \cite{Bilokopytov:21}. They are related to the tauberian operators, defined by Kalton and Wilansky \cite{KaltonWil:76}; in fact, they coincide when $E=L_1$ (see \cite{GM:97} and \cite{GMM:20}). We refer to \cite{GO:90} and \cite{GM:10} for additional information on tauberian operators. 

The disjointly non-singular operators can be characterized as follows.

\begin{Thm}\label{thm:Lp-D} \cite[Theorem 2.8]{GMM:20}
For $T\in\Lc(E,Y)$, the following assertions are equivalent:
\begin{enumerate}
\item $T$ is disjointly non-singular.
\item There is no disjoint sequence of non-zero vectors $(x_n)$ in $E$ such that the restriction $T_{[x_n]}$ is a  compact operator.

\item For every disjoint sequence of non-zero vectors $(x_n)$ in $E$, the restriction $T_{[x_n]}$ is an upper semi-Fredholm operator.

\item For every normalized disjoint sequence $(x_n)$ in $E$, $\liminf_{n\rightarrow \infty} \|Tx_n\|>0$.

\end{enumerate}
\end{Thm}

It was proved in \cite[Proposition 14]{GM:97} and \cite[Theorem 3.15]{GMM:20} that, for $1\leq p<\infty$, $\dss(L_p,Y)$ is the perturbation class of $\dns(L_p,Y)$.

\subsection*{Representation of Banach lattices}\label{rep-BL}
It is well-known (see \cite[Theorem 1.b.14]{LT-II:79}) that every order continuous Banach lattice with a weak unit $E$ admits a representation as a K\"othe function space, in the sense that there exists a probability space $(\Omega,\Sigma,\mu)$ so that
\begin{itemize}
\item $L_\infty(\mu) \subset E \subset L_1(\mu)$ with $E$ dense in $L_1(\mu)$ and $L_\infty(\mu)$ dense in $E$, 
\item $\|f\|_1\leq \|f\|_E\leq 2\|f\|_\infty$ when $f\in L_\infty(\mu)$, 
\end{itemize}
and the order in $E$ is the order induced by $L_1(\mu)$. 

The following fact will allow us to state some of our results omitting the existence of a weak unit in the Banach lattice. 

\begin{Lemma}\label{no-unit}  
Let  $E$ be an order continuous Banach lattice. Then each sequence in $E$ is contained in a closed ideal of $E$ with a weak unit. 
\end{Lemma} 
\begin{proof}
If $(f_n)$ is a bounded sequence in $E$, then $e=\sum_{n=1}^\infty |f_n|/2^{-n}$ is a weak unit in the closed ideal generated by $(f_n)$. 
\end{proof}

We also will need the following result.  

\begin{Lemma}\label{rem:ord-cont}  Let $E$ be an  order continuous Banach lattice with a weak unit, and let $f \in E$. If $(A_k)$ is a disjoint sequence in the $\sigma$-algebra $\Sigma$ associated to the representation of $E$, then $\lim_{k\to\infty}\|f 1_{A_k}\|_E=0$.  
\end{Lemma}
\begin{proof}
Let $B_k=\cup_{i=k}^\infty A_i$. Since the norm on $E$ is order continuous, $(B_k)$ is decreasing and $\lim_{k \rightarrow \infty} \mu (B_k) =0$ we have $\lim_{k \to \infty} \|f 1_{B_k} \|_E = 0$, hence $\lim_{k\to\infty} \|f 1_{A_k}\|_E= 0$.
\end{proof}


\section{Operational quantities} 

An \emph{operational quantity} is a map $a:\Lc(X,Y)\to[0,\infty)$ satisfying certain conditions. Given two operational quantities $a$ and $b$, we write $a\leq b$ when $a(T)\leq b(T)$ for each $T\in \Lc(X,Y)$. Moreover, the quantities $a$ and $b$  are {\it equivalent} if there exist positive constants $c_1<c_2$ such that $c_1 a \leq b \leq c_2 a$.

We are interested in some classical operational quantities and some new ones that we introduce here. To describe the classical ones, let $S(X)$ be set of all closed infinite dimensional subspaces of $X$. Then, given an operational quantity $a:\Lc(X,Y)\to[0,\infty)$, we define two derived quantities $i\, a$ and $s\, a$ as follows:
\begin{equation}\label{ia,sa}
i\, a(T) :=\inf_{M\in S(X)} a(T_M)\;\; \mbox{ and } \;\; s\, a(T) := \sup_{M \in S(X)} a(T_M), 
\end{equation}
 where $T\in\Lc(X,Y)$.

Note that $a \leq b$ implies $i a \leq i b$ and $s a \leq s b$. Taking the operator norm as $a$ in (\ref{ia,sa}), for $T\in\Lc(X,Y)$ we obtain  
\begin{itemize}
  \item $\G(T) := i\, \|T\| = \inf_{M\in S(X)} \|T_M\|$ \; \;  and \smallskip
\item $\D(T) := s\, \G(T) = \sup_{M\in S(X)} \G(T_M)= \sup_{M\in S(X)} \inf_{N\in S(M)} \|T_N\|$.\smallskip
\end{itemize}
The quantities $\G = i\, \| \cdot \|$ and $\D = i\, \G$ were introduced by Gramsch and Schechter (see \cite{schechter,schechter book}), who proved that $\G(T)>0$ if and only if $T$ is upper semi-Fredholm, and $\D(T)=0$ if and only if $T$ is strictly singular.

To introduce the new quantities, we denote by $\ds(E)$ the set of all sequences of disjoint non-zero vectors of $E$. Now, given an operational quantity $a:\Lc(F,Y) \to[0,\infty)$ defined for $F=E$ and $F\in \ds(E)$, for each $T\in\Lc(E,Y)$ we define two derived quantities $i_d\, a$ and $s_d\, a$ as follows: 
\begin{equation}\label{ida, sda} 
i_d\, a(T) := \inf_{(x_n) \in \ds(E)} a(T_{[x_n]})\quad \mbox{ and } \quad s_d\, a(T) := \sup_{(x_n) \in \ds(E)} a(T_{[x_n]})  \, . 
\end{equation}
Again, $a \leq b$ implies $i_d a \leq i_d b$ and $s_d a \leq s_d b$.
We are interested in two operational quantities derived from the norm, whose notation is inspired by that of Schechter:
\begin{itemize}
  \item $\G_d(T) := i_d\, \|T\| = \inf_{(x_n)\in \ds(E)} \|T_{[x_n]}\|$\quad and\smallskip
  \item $\D_d(T) :=s_d\, \G_d(T)= \sup_{(x_n)\in \ds(E)} \G_d (T_{[x_n]}) = \sup_{(x_n)\in \ds(E)}\inf_{(y_n) \in\ds([x_n])} \|T_{[y_n]}\|$,\smallskip 
\end{itemize}
that will allow us to characterize the operators in $\dns$ and $\dss$.

In a similar way, for $T\in\Lc(X,Y)$ we consider two  classical operational quantities derived from the injection modulus $j$: 

\begin{itemize}
  \item $\tau(T) := s\, j(T) = \sup_{M\in S(X)} j(T_M)$ \; \;  and\smallskip 
\item $\s(T) :=i\, \tau(T)= \inf_{M\in S(X)} \tau(T_M)  =\inf_{M\in S(X)} \sup_{N\in S(M)} j(T_N) $, 
\end{itemize}
and derive two new quantities for $T\in\Lc(E,Y)$: 
\begin{itemize}
  \item $\tau_d(T) := s_d\, j(T) = \sup_{(x_n) \in \ds (E)} j(T_{[x_n]})$ \; \;  and \smallskip
\item $\s_d(T) := i_d\, \tau_d(T)= \inf_{ (x_n) \in \ds (E) } \tau_d(T_{[x_n]})  = \inf_{(x_n) \in \ds (E)} \sup_{(y_n) \in \ds ([x_n])} j(T_{[y_n]})$,\smallskip 
\end{itemize}

The operational quantities $\tau = s\, j$ and $\s = i\, \tau$ were introduced in \cite{schechter} and \cite{GM:93}, where it was proved that $\tau(T)=0$ if and only if $T$ is strictly singular, and $\s(T)>0$ if and only if $T$ is upper semi-Fredholm. 
We will show that the quantities  $\tau_d$ and $\s_d$ characterize the operators in $\dss$ and $\dns$, respectively. 
\medskip

The proof of the next lemma shows that for each closed infinite dimensional subspace of a Banach space with a monotone basis $(x_n)$, in particular with a  $1$-unconditional basis, there is a block basis $(y_k)$ such that $[y_k]$ is `arbitrarily close' (in the sense of the gap between subspaces; see  \cite[Section IV.2]{Kato:80}) to a subspace $N$ of $M$; so the action of an operator on $[y_k]$ is also close to its action on $N$. This idea will appear several times in our arguments. 

\begin{Lemma}\label{lemma} 
Let $X$ be a Banach space with a monotone basis $(x_n)$, let $M\in S(X)$ and $0<\e<1$. Then there exist a normalized block basis $(y_k)$ of $(x_n)$ and a subspace $N\in S(M)$ such that for every operator $T \in \Lc(X,Y)$, 
$$
\left|\| T_{[y_k]}\| -\| T_N\| \right| \leq \e\|T\| \;\; \mbox{ and }\;\; \left|j(T_{[y_k]})- j(T_N) \right|\leq \e\|T\|.
$$
\end{Lemma}
\begin{proof}  
We will choose $(y_k)$ and $N$ so that the distance between the unit spheres of $N$ and $[y_k]$ is smaller than $\e$; hence for each $n\in S_N$ there is $y\in S_{[y_k]}$ with $\|n-y\|< \e$, and for each $z\in S_{[y_k]}$ there is $m\in S_N$ with $\|z-m\|< \e$. Clearly this fact implies our result.

Let $r=\e/8$. Inductively, we will find integers $1=j_1\leq l_1<j_2\leq l_2 \leq\cdots$  and a sequence $(a_i)$ of scalars so that $y_k= \sum_{i=j_k}^{l_k} a_i x_i$ satisfies $\|y_k\|=1$ and $\dist(y_k,M)<r/2^{k+1}$. 

Clearly, $y_1$ exists; so assume that $y_k$ has been found for $k\leq k_0$. Let $(x^*_i)$ be the sequence in $X^*$ such that $x^*_i(x_j)= \delta_{i,j}$. Since $M\cap \left(\cap_{i=1}^{l_{k_0}} N(x^*_i)\right)$ is infinite dimensional, $y_{k_0+1}$ exists. 

Since $(y_k)$ is a monotone basic sequence (comment after \cite[Definition 1.a.10]{LT-I:77}),  there exists a sequence $(y_k^*)$ in  $X^*$ with $\|y_k^*\|\leq2$ and $y_k^*(y_j)=  \delta_{k,j}$.

For each $k\in \N$ we choose $m_k\in M$ with  $\|y_k-m_k\|< r/2^{k+1}$, and define $K\in \Lc(X)$ by 
$$
Kx := \sum_{k=1}^\infty y_k^*(x) (y_k-m_k).
$$
Then $K$ is bounded with $\|K\| \leq \sum_{k=1}^\infty \|y^*_k\|\cdot\|y_k-m_k\|< r$; hence $I-K$ is bijective. Moreover $(I-K)y_k=m_k$ for each $k\in\N$. We take $N=[m_k]=(I-K)([y_k])$. 
Note that 
$$
(I-K)^{-1} = \sum_{l=0}^\infty K^l = I-L \textrm{ with }\|L\|\leq \sum_{l=1}^\infty r^l = r/(1-r)<2r. 
$$

For $n\in S_N$ we take $y= \|(I-L)n\|^{-1}(I-L)n\in S_{[y_k]}$. Then $1-2r <\|(I-L)n\|<1+2r$ and  
$$
\|n-y\|=\frac{\big\|(\|(I-L)n\|-1)n+Ln\big\|}{\|(I-L)n\|} \leq \frac{4r}{1-2r}<8r=\e.
$$

Similarly, for each $z\in S_{[y_k]}$, we have $m=\|(I-K)z\|^{-1}(I-K)z\in S_N$ and  $\|z-m\|< \e$.
\end{proof}

A Banach lattice is called  \emph{atomic} if its order is induced by a $1$-unconditional basis.

\begin{Prop}\label{E ub} Let $E$ be an atomic Banach lattice. For an operator $T\in\Lc(E,Y)$, 
$$
\G_d (T) = \G(T)  \; ,  \;  \Delta_d (T) = \Delta (T)  \; , \; \tau_d (T) = \tau(T)  \; \mbox{ and }  \;  \s_d (T) = \s (T)  \; .  
$$ 
\end{Prop}
\begin{proof} The inequality  $\G_d (T) \geq \G(T)$ is valid in general. The converse inequality is obtained by applying Lemma \ref{lemma}. Suppose without loss generality that $\| T \| =1$. Given $0 < \e < 1$ and a subspace $M$ of $E$, there is a block basis $(y_k)$ of the unconditional basis of $E$ such that $[y_k]$ is arbitrarily close to some subspace $N$ of $M$, and consequently  
$$
\left| \| T_{[y_k]} \| - \| T_{N} \| \right| \leq \e  \; . 
$$
Hence $\G_d(T)\leq \|T_{[y_k]} \| \leq \| T_{N} \| + \e   \leq \| T_{M} \| + \e$.
Therefore $\G_d(T) \leq \G(T)$.

The other equalities can be  proved in a similar way. 
\end{proof} 

\begin{Cor}\label{cor:i_ds_d=i_ds}
We have $s_d\, \G_d=s_d\, \G$ and $i_d\, \tau_d = i_d\, \tau$. Moreover $\G_d=i_d\, \G_d=i_d\, \G$ and $\tau_d= s_d\, \tau_d= s_d\, \tau$. 
\end{Cor}
\begin{proof}
For each $(x_n)\in \ds (E)$, $(x_n)$ is a $1$-unconditional basis; hence $[x_n]$ is an atomic Banach lattice. Therefore 
$$s_d\, \G_d(T)= \sup_{(x_n)\in \ds(E)} \G_d (T_{[x_n]}) = \sup_{(x_n)\in \ds(E)} \G (T_{[x_n]})= s_d\, \G(T).$$

The proof of $i_d\, \tau_d = i_d\, \tau$, $i_d\, \G_d=i_d\, \G$ and $s_d\, \tau_d= s_d\, \tau$  is identical, and for the remaining equalities, note that $i_d\, i_d\, a =i_d\, a$ and $s_d\, s_d\, a =s_d\, a$ for any quantity $a$. 
\end{proof}


\section{Operational quantities derived from the norm}

Our first result gives some alternative expressions for $\G_d (T)$ in terms of the classical quantities.

\begin{Prop}\label{prop:i_D n} For $T\in\Lc(E,Y)$, we have
$\G_d (T) = i_d \G (T) = i_d \D (T)$.
\end{Prop}
\begin{proof} 
Note that $\G_d= i_d\, \|\cdot \|$. Applying $i_d$ to the inequalities $\G \leq \D \leq \| \cdot \|$, we obtain $i_d\, \G\leq i_d\, \D\leq i_d\, \|\cdot \| $, and Corollary \ref{cor:i_ds_d=i_ds} completes the proof.
\end{proof} 

It was proved in \cite{GMM:20} that $T\in\Lc(E,Y)$ is disjointly non-singular if and only if for every $(f_n)\in\ds(L_p)$, the restriction $T_{[f_n]}$ is upper semi-Fredholm.
Next we give a quantitative version of this result when $E$ is an order continuous Banach lattice.
Since $\G_d (T) = i_d \G (T)$ by Proposition \ref{prop:i_D n}, our result says that if $T\in\dns(E,Y)$ then the restrictions $T_{[x_n]}$ are ``uniformly" upper semi-Fredholm, in the sense that  $\inf_{(x_n)\in \ds(E)} \G(T_{[x_n]}) > 0$.  

\begin{Thm}\label{th:gamma-d}
Let $E$ be an order continuous Banach lattice, and let $T \in L(E,Y)$. Then $T\in \dns$ if and only if $\G_d (T)>0$.
\end{Thm}

\begin{proof} 
Suppose that $\G_d (T) >0$. For every $(f_n)\in\ds(E)$ we have that $\G(T_{[f_n]}) > 0$,
hence $T_{[f_n]}$ is upper semi-Fredholm.
Consequently, $T$ is disjointly non-singular (Theorem \ref{thm:Lp-D}).

Conversely, we assume that $\G_d (T)=0$. By Theorem \ref{thm:Lp-D}, it is enough to construct a normalized sequence $(h_n)\in \ds(E)$ such that $\lim_{n\to\infty}\|Th_n\|=0$. 

For each $n\in\N$ there exists  a normalized sequence $(f_{n,k})_k\in\ds(E)$ such that
$\|T_{[(f_{n,k})_k]} \| < 1/n$, and by Lemma \ref{no-unit} we can assume that the functions $f_{n,k}$ ($n,k\in\N$) are contained in a closed ideal of $E$ which has a representation as a K\"othe space. 

Let $g_1 = f_{1,1}$. As $\lim_{k\to\infty} \mu (D(f_{2,k})) = 0$, by Lemma  \ref{rem:ord-cont} we have  $\lim_{k\to\infty}\|g_1 1_{D(f_{2,k})} \|_E =0$. So we can find  $k_2>1$ such that
$$
\| g_1 \| = 1 \; , \; \|Tg_1 \| < 1  \; \mbox{ and } \; \|g_1 1_{D(f_{2,k_2})} \|_E < \frac{1}{2^2}. 
$$

Then, taking $g_2 = f_{2,k_2}$, a similar argument using Lemma \ref{rem:ord-cont} shows that there exists $k_3>k_2$ such that 
$$
\| g_2 \| = 1 \; , \; \|Tg_2 \| < \frac{1}{2}  \; \mbox{ and } \;\|g_i 1_{D(f_{3,k_3})} \|_E < \frac{1}{2^3}\quad \textrm{ for $1 \leq i < 3$.} 
$$

In this way we find a sequence $k_1 = 1 < k_2 < k_3 < \cdots$ such that, taking $g_l = f_{l,k_l}$ for each $l\in\N$, we have
$$
\| g_l\| = 1 \; , \; \|Tg_l\| < \frac{1}{l}  \; \mbox{ and } \; \|g_i 1_{D(f_{l,k_{l+1}})} \| < \frac{1}{2^{l+1}} \; \; \; ( 1\leq i < l+1).
$$

Let $A_k = \cup_{j=k+1}^\infty D(g_j)$ and $\tilde{h}_k := g_k - g_k 1_{A_k}$. For $k < l$ we have $D(\tilde{h}_k) \cap D(g_l) = \emptyset$ and $D(\tilde{h}_l ) \subset D(g_l)$, hence $D(\tilde{h}_k) \cap D(\tilde{h}_l) =  \emptyset$. Thus the sequence $(\tilde{h}_k)$ is disjoint.
Since $\| g_n \| = 1$,
\begin{eqnarray*}
|1 - \| \tilde{h}_n \| | &\leq& 
\|g_n- \tilde{h}_n\|= \| g_n 1_{A_n} \|\\ 
&\leq& \left\| \sum_{i=n+1}^\infty  g_n 1_{D(g_i)} \right\| \leq \sum_{i=n+1}^\infty \| g_n 1_{D(g_i)}\|\\ 
&\leq& \sum_{i=n+1}^\infty \frac{1}{2^i} = \frac{1}{2^n}.
\end{eqnarray*}

Taking $h_n= \|\tilde{h}_n\|^{-1} \tilde{h}_n$, we obtain $(h_n)\in \ds(E)$ is normalized and  
\begin{eqnarray*}
\|h_n - g_n \| &\leq& 
\left\|\frac{\tilde{h}_n}{\|\tilde{h}_n\|} - \frac{g_n}{\| \tilde{h}_n \|} \right\|
+ \left\| \frac{g_n}{\| \tilde{h}_n \|} - g_n \right\|\\
&=& \frac{\| \tilde{h}_n - g_n \|}{\| \tilde{h}_n \|} + \frac{ |1-\| \tilde{h}_n\| | \; \|g_n\|}{\| \tilde{h}_n \|}\\
&\leq&
\frac{2\|\tilde{h}_n-g_n\|}{\| \tilde{h}_n\|}
\leq \frac{1}{2^{n-1} \|\tilde{h}_n\|}.
\end{eqnarray*}
Consequently $\lim_{n\to\infty} \|h_n- g_n\|=0$, and $\|Th_n\|\leq \|T(h_n - g_n)\| + \|Tg_n\|$ and $\|Tg_n\|<1/n$; hence $\lim_{n\to\infty} \|T h_n\|= 0$.
\end{proof}

Next we give some alternative expressions for $\D_d (T)$.

\begin{Prop}\label{prop:s_D n}
For $T\in\Lc(E,Y)$, we have $\D_d (T) = s_d \D(T) = s_d \G (T)$.
\end{Prop}
\begin{proof} 
Note that $\D_d (T) = s_d \G_d (T)$ and, by Corollary \ref{cor:i_ds_d=i_ds}, $s_d\, \G (T) = s_d\, \G_d (T)$. So it is enough to observe that $s_d\, a (T)= s_d\, s_d\, a$ for any quantity $a$.
\end{proof} 

\begin{Prop}\label{prop:DSS-n}  $T\in\Lc(E,Y)$ is disjointly strictly singular if and only if $\D_d (T) =0$.
\end{Prop}
\begin{proof}
As $\D_d (T) = s_d \D (T)$, we have that $\D_d (T) = 0$ means that for every $(x_n) \in \ds(E)$ we have that $\D(T_{[x_n]}) =0$; that is, all the restrictions  $T_{[x_n]}$ are strictly singular. By \cite[Proposition 2.6]{GMM:20}, that is equivalent to $T$ being disjointly strictly singular.
\end{proof}

Obviously, given $T\in\Lc(E,Y)$ and a scalar $\lambda$, $\G_d(\lambda T)=|\lambda| \G_d(T)$ and
$\D_d(\lambda S)=|\lambda|\D_d(S)$.
The following result complements these facts.

\begin{Prop}\label{prop:ineqn2} 
For operators $T,S\in\Lc(E,Y)$, we have the following inequalities:
\begin{enumerate}
\item $\G_d (T+S)\leq \G_d (T)+ \D_d  (S)$ and
\item $\D_d (T+S) \leq \D_d ( T) + \D_d (S)$.
\end{enumerate}
\end{Prop}

\begin{proof} Let $(x_n)\in\ds (E)$. Then  $\|(T+S)_{[x_n]}\|\leq \|T\|+ \| S_{[x_n]} \|$, and  taking the infimum over $(x_n)\in\ds (E)$ we obtain $\G_d (T+S)\leq \|T\| + \G_d (S)$. Therefore 
$$
\G_d (T+S)\leq\G_d \left((T+S)_{[x_n]} \right)\leq \|T_{[x_n]}\|+ \G_d (S_{[x_n]}) \leq \|T_{[x_n]}\| + \D_d (S),
$$
and taking again the infimum over $(x_n)\in\ds (E)$ we get (1).

Let $(x_n)\in\ds (E)$. From (1) we derive $$
\G_d ((T+S)_{[x_n]}) \leq \G_d (T_{[x_n]})+ \D_d(S_{[x_n]})\leq \G_d(T_{[x_n]})+ \D_d(S),
$$ 
and taking the supremum over $(x_n)$ we get $\D_d (T+S) \leq \D_d (T) + \D_d (S)$.
\end{proof}

Since $\D_d(T)\leq\|T\|$, Theorem \ref{th:gamma-d} and part (1) of Proposition \ref{prop:ineqn2} improve the results proved in \cite{GMM:20} that, under some conditions,  $\dns(E,Y)$ is stable under perturbation by small norm operators and $\dss$ operators. 

\begin{Cor}\label{prop:pert n} 
Let $E$ be an order continuous Banach lattice. Then  
\begin{enumerate}
\item $\dss (E,Y)$ is a closed subspace of $L(E,Y)$;
\item $\dns (E,Y)$ is an open subset of $L(E,Y)$;
\item If $S\in \dss(E ,Y)$, then $\G_d(T+S) =\G_d(T)$, for all $T\in\Lc(E,Y)$;\newline
in particular, $T \in \dns (E,Y)$ implies $T+S \in \dns(E,Y)$.
\end{enumerate}
\end{Cor}
\begin{proof} (1) If $T,S \in \dss(E,Y)$, then $\D_d (T+S) \leq \D_d (T) + \D_d (S)=0$, so $T+S \in \dss(E,Y)$; and $\D_d (\la T) =|\la|\D_d (T)$ implies $\la T\in  \dss(E,Y)$.

(2) If $T \in \dns (E,Y)$ and $S \in \Lc(E,Y)$ with $\|S\| < \G_d(T)$, then $\G_d(T+S) \geq \G_d (T) - \D_d (S)\geq \G_d(T)- \|S\|>0$. Hence $T+S\in \dns (E,Y)$. 

(3) Let $S \in \dss(E,Y)$, so $\D_d(S)=0$. For all $T \in \Lc(E,Y)$, 
$$
\G_d (T+S) \leq \G_d(T) + \D_d (S) = \G_d (T), 
$$
and similarly $\G_d(T)= \G_d(T+S-S)\leq \G_d(T+S)$. 
\end{proof} 

Part (2) of Corollary \ref{prop:pert n} was proved by Bilokopytov \cite{Bilokopytov:21} using different techniques.
\medskip

A closed subspace $M$ of $E$ is said to be {\it dispersed} if there is no sequence $(x_n)\in\ds(E)$ such that $\lim_{n\to\infty}\dist(x_n,M)=0$ (see \cite[Definition 2.1]{GMM:20}). 

\begin{Rem}\label{rem:dispersed-n} 
Let $M$ be a non-dispersed closed subspace of $E$. Denoting by $\textrm{ND}(M)$ the set of all closed subspaces of $M$ which are non-dispersed in $E$, it readily follows from Lemma \ref{lemma} that, for $T \in \Lc(E,Y)$,  
$$
\G_d(T) = \inf_{M \in \textrm{ND}(E)} \|T_M\| \quad \mbox{ and } \quad \D_d(T) = \sup_{M_1 \in \textrm{ND}(E)}\; \inf_{M_2 \in \textrm{ND}(M_1)}\|T_{M_2}\|. 
$$
\end{Rem}


\section{Operational quantities derived from the injection modulus}

Next result gives other expressions for the quantity $\tau_d$. 

\begin{Prop}\label{prop:s_Dj} For $T \in \Lc(E,Y)$, we have $\tau_d (T)= s_d \s (T) = s_d \tau (T)$. 
\end{Prop}
\begin{proof} As $j \leq \s \leq \tau$, we have $\tau_d = s_d j \leq s_d \s \leq s_d \tau$. Moreover, $s_d\, \tau =s_d\, \tau_d$ by Corollary \ref{cor:i_ds_d=i_ds}. Hence 
$$s_d\, \tau (T)=s_d\, \tau_d(T)=s_d\, s_d\, j(T)= s_d\, j(T)= \tau_d(T),$$ 
because $s_d\, s_d\, a= s_d\, a$ for every quantity $a$. 
\end{proof} 

\begin{Prop}\label{prop:DSS-j} 
Let $T\in\Lc(E,Y)$. Then $T\in\dss$ if and only if $\tau_d(T) =0$.
\end{Prop}

\begin{proof}
We have that $\tau_d (T) = 0$ is equivalent to $j(T_{[x_n]}) = 0$, for every sequence $(x_n) \in \ds(E)$. This means that $T$ is not an isomorphism on any subspace $[x_n]$ generated by a disjoint sequence. That is, $T$ is  disjointly strictly singular.
\end{proof}

\begin{Prop}\label{prop:i_Dj} For an operator $T\in\Lc(E,Y)$, we have 
$\s_d (T) = i_d \s (T) = i_d \tau (T)$. 
\end{Prop}
\begin{proof} By Proposition \ref{prop:s_Dj},  $\s \leq \tau_d \leq \tau$, hence $i_d \s \leq i_d \tau_d = \s_d \leq i_d \tau$. Moreover, arguing as in the proof of Corollary \ref{cor:i_ds_d=i_ds} we get $i_d \s =i_d\,  \s_d=i_d\, i_d\, \tau_d= i_d\, \tau_d= i_d\, \tau$, and the result is proved. 
\end{proof}

Like Theorem \ref{th:gamma-d}, by Proposition \ref{prop:i_Dj} the following result says that $T\in\dns(E,Y)$ if and only if the restrictions $T_{[x_n]}$ with $(x_n)\in \ds(E)$ are ``uniformly" upper semi-Fredholm, in the sense that $\inf_{(x_n)\in\ds(E)} \s(T_{[x_n]})> 0$. 

\begin{Thm}\label{thm:DNS j} Let $E$ be an order continuous Banach lattice and let $T \in L(E,Y)$. Then $T\in \dns$ if and only if $\s_d (T)>0$.
\end{Thm}

\begin{proof} 
By Proposition \ref{prop:i_Dj}, $\s_d (T) = i_d \tau (T)$. Then if $\s_d(T)>0$ and $(f_n)\in \ds(E)$, $\tau (T_{[f_n]}) > 0$. Hence $T_{[f_n]}$ is not  strictly singular, and $T$ is disjointly non-singular by Theorem \ref{thm:Lp-D}.

Conversely, suppose that $\s_d (T)=0$. By Theorem \ref{thm:Lp-D}, in order to  show that $T$ is not disjointly non-singular, it is enough to find a normalized $(h_n) \in \ds(E)$ such that $\lim_{n\to\infty} Th_n= 0$. 

For each $n\in\N$ there exists a normalized sequence $(f_{n,k})_k \in \ds (E)$  such that
$$
\tau_d (T_{[f_{n,k}]_k} ) < \frac{1}{n},
$$
and by Lemma \ref{no-unit} we can assume that the vectors $f_{n,k}$ are contained in a closed ideal that admits a representation as a K\"othe space.

As  $j(T_{[f_{1,k}]_k} ) < 1$, there exists $g_1 \in [(f_{1,k})_k]$  with $\| Tg_1 \| < 1$. From $\lim_{k \to \infty} \mu (D(f_{2,k})) = 0$, by Lemma  \ref{rem:ord-cont} we have  $\lim_{k\to\infty}\|g_1 1_{D(f_{2,k})} \|_E =0$. So we can to take $k_2>1$ such that
$$
\| g_1 \| = 1 \; , \; \|Tg_1 \| < 1  \; \mbox{ and } \; \|g_1 1_{D(f_{2,k_2})} \|_E < \frac{1}{2^2}. 
$$
Moreover, from
$$
j(T_{[(f_{2,k})_{k \geq k_2}]} ) \leq \tau_d (T_{[(f_{2,k})_k]} ) < \frac{1}{2} \; ,
$$
we obtain that there is $g_2 \in [(f_{2,k})_{k \geq k_2}]$ with $\| Tg_2 \| < 1/2$. As $\lim_{k \to \infty} \mu (D(f_{3,k})) = 0$, by Lemma \ref{rem:ord-cont} we get   $\lim_{k\to\infty}\|g_i 1_{D(f_{3,k})} \|_E =0$, so we can take $k_3>k_2$ such that
$$
\| g_2 \| = 1 \; , \; \|Tg_2 \| < \frac{1}{2}  \; \mbox{ and } \; \|g_i 1_{D(f_{3,k_3})} \|_E < \frac{1}{2^3} \; (i \leq i < 3) \; .  
$$

Now, proceeding as in the proof of Theorem  \ref{th:gamma-d}, we take $A_n = \cup_{j=n+1}^\infty D(g_j)$ and obtain a normalized sequence $h_n := \|g_n-g_n 1_{A_n} \|^{-1}(g_n - g_n 1_{A_n})$ in $\ds(E)$. Since $\lim_{n\to\infty}\|T h_n\|= 0$, we conclude that $T\notin \dns(E,Y)$. 
\end{proof}

To compare Theorem \ref{thm:DNS j} with  Theorem \ref{th:gamma-d}, observe that $\s_d \leq \G_d$.

\begin{Prop}\label{lemma ineq j} For operators $T,S\in\Lc(E,Y)$, we have the following inequalities:
\begin{enumerate}
\item $\tau_d (T+S) \leq \tau_d ( T) + \D_d (S)$ and 
\item $\s_d (T+S) \leq \s_d ( T) + \D_d(S)$. 
\end{enumerate}
\end{Prop}

\begin{proof} Since $j(T+S)\leq j(T)+\|S\|$,  for each $(x_n) \in \ds(E)$ we get 
$$j(T+S)\leq j((T+S)_{[x_n]}) \leq j(T_{[x_n]}) + \|S_{[x_n]}\|\leq \tau_d (T)+ \|S_{[x_n]}\|,$$ 
and taking the infimum over $(x_n)$ we obtain  $j(T+S) \leq \tau_d (T) + \G_d (S)$.

(1) For $(x_n) \in \ds(E)$, we have $j((T+S)_{[x_n]}) \leq \tau_d(T_{[x_n]}) + \G_d (S_{[x_n]})\leq \tau_d(T) + \G_d (S_{[x_n]})$, and taking the supremum over  $(x_n)$ we get $\tau_d(T+S) \leq \tau_d(T)+ \D_d (S)$.

(2) Applying (1), $\tau_d ((T+S)_{[x_n]}) \leq \tau_d(T_{[x_n]})+ \D_d(S_{[x_n]})\leq \tau_d ( T_{[x_n]}) + \D_d(S)$ for each  $(x_n) \in \ds(E)$. So taking the infimum over $(x_n)$, we obtain $\s_d(T+S)\leq \s_d(T)+ \D_d(S)$.
\end{proof}

From Proposition \ref{lemma ineq j}, we could derive an alternative proof of Corollary \ref{prop:pert n}.

\begin{Rem}\label{rem:dispersed-j} 
As in Remark \ref{rem:dispersed-n}, we can give expressions for $\s_d(T)$ and $\tau_d(T)$ in terms of the restrictions of $T$ to non-dispersed subspaces. 
For $T \in \Lc(E,Y)$,

$\tau_d(T)= \sup_{M\in \textrm{ND}(E)} j(T_M)$ \quad and  \quad $\s_d(T) = \inf_{M_1 \in \textrm{ND}(E)} \sup_{M_2 \in \textrm{ND}(M_1)} j(T_{M_2})$. 
\end{Rem}


\section{The quantity $\beta$}

For an operator  $T\in\Lc(E,Y)$, the following quantity was defined in \cite{GMM:20}:
$$
\beta(T) :=\inf\left\{\liminf_{n\to\infty} \|Tx_n\| : \textrm{ $(x_n)$ normalized disjoint in }E\right\}.
$$


We have shown in Theorem \ref{th:gamma-d} that the quantity $\G_d$ characterizes $\dns(E,Y)$ for $E$ an order continuous Banach lattice. Moreover, it is related with $\beta$ as follows:
\begin{Prop}\label{beta1}
Every operator $T \in L(E,Y)$ satisfies $\beta(T) \leq \G_d (T)$.
\end{Prop}
%
\begin{proof}
Note that 
\newline\indent
$\beta(T) = \inf_{(x_n)\in \ds(E)}\liminf_{n\to  \infty} \left\|T\frac{x_n}{\|x_n\|}\right\|
\leq\inf_{(x_n)\in\ds(E)} \|T_{[x_n]}\| =\G_d (T)$.
\end{proof}
%
%
It was proved in \cite[Proposition 3.1]{GMM:20}  (see \cite{GM:97} for $p=1$) that, for $1\leq p<\infty$, an operator $T \in \Lc(L_p,Y)$ is disjointly non-singular if and only if $\beta (T)>0$. Now we extend this result. 

%

\begin{Prop}\label{beta3}
Let $E$ be an order continuous Banach lattice. Then an operator  $T \in L(E,Y)$ is disjointly non-singular if and only if $\beta (T)>0$.
\end{Prop}
\begin{proof} 
If $\bb(T)>0$, then condition (4) in Theorem \ref{thm:Lp-D} is satisfied, hence $T \in \dns (E,Y)$. 
\newline\indent
Suppose that $\bb(T)=0$. Then for every $n\in\N$ we can find a normalized disjoint sequence $(f_{n,k})_{k\in\N}$ with $\|Tf_{n,k}\|<1/n$ for every $k\in \N$, and proceeding as in the proof of Theorem \ref{th:gamma-d}, for each $n$ we select $k_n$ so that taking $g_n= f_{n,k_n}$ we have $\|g_i 1_{D(g_n)}\|<2 ^{-n}$ for $1\leq i <n$. 
The sequence $(g_n)$ is almost disjoint (there exists a normalized disjoint sequence $(h_n)$ in $E$ such that $\lim_{n\to\infty} \|g_n-h_n\|_E=0$).   Then $\lim_{n\to\infty}\|Th_n\|=0$, hence $T \notin \dns (E,Y)$.  
\end{proof}

By Proposition \ref{beta1}, $\beta \leq \G_d$. In some cases, these two  quantities coincide; for example, if $1\leq p<2$ and $M$ is a dispersed subspace of $L_p$, then the quotient map $Q_M:L_p\to L_p/M$ satisfies $\beta(Q_M)=1$ (see \cite{GMM:20}), hence $\G_d(Q_M)=\|Q_M\|=1$.
However, using the fact proved by Odell and Schlumprecht in  \cite{OdellS:94} that the Banach space $\ell_2$ is arbitrarily distortable, we show that these two quantities are not equivalent:

\begin{Example}\label{ex:distor}  For every $\la>1$ and $\e>0$, there exists a Banach space $Y_\la$ isomorphic to $\ell_2$ and an operator $T_\la\in\Lc(\ell_2,Y_\la)$ such that $0<\la\cdot\beta(T_\la)\leq \G_d(T_\la)+\e$.  Thus there is no $C>0$ such that $\G_d\leq C\cdot \beta$. 
\end{Example}
\begin{proof}
Since $\ell_2$ is arbitrarily distortable \cite{OdellS:94}, for every $\la>1$ there is a norm $|\cdot|_\la$ on $\ell_2$ equivalent to the usual one $\|\cdot\|_2$ such that, for each closed infinite dimensional subspace $M$ of $\ell_2$, 
\begin{equation}\label{eq:indesc}
\sup\left\{\frac{|x|_\la}{|y|_\la} : x,y\in M, \|x\|_2=\|y\|_2=1\right\} > \la. 
\end{equation}
We denote  $Y_\la=(\ell_2, |\cdot|_\la)$ and $T_\la$ the identity operator from $\ell_2$ onto $Y_\la$. 

Note that the operator $T_\la$ is bounded below, and passing to a closed infinite dimensional subspace of $\ell_2$ (that we can identify with $\ell_2$, with the lattice structure determined by any orthonormal basis) we can assume that $\|T_\la\| < \G_d(T_\la)+\e$. 

By inequality (\ref{eq:indesc}), $\la j(T_\la)\leq \|T_\la\|$ and there exists $g_1$ with $\|g_1\|_2=1$ and $\la\cdot |g_1|_\la <\G_d(T_\la) +\e$. Moreover, by the denseness of the span of the basis $(e_n)$ of  $\ell_2$, we can choose $g_1\in [e_1,\ldots, e_{m_1}]$ for some $m_1\in\N$. 
Similarly, there exists $g_2\in [e_i :i>m_1]$ with $\|g_2\|_2=1$ and $\la\cdot |g_2|_\la <\G_d(T_\la) +\e$, and again we can choose $g_2\in [e_{m_1+1},\ldots, e_{m_2}]$ for some $m_2>m_1$ in $\N$. 

In this way we get a sequence $(g_n)\in \ds(\ell_2)$ such that $\la\cdot|g_n|_\la= \la\cdot|T_\la g_n|_\la\leq \G_d(T_\la)+\e$, which implies  $\la\cdot \beta(T_\la)\leq \G_d(T_\la) +\e$.
\end{proof}



\section{Order between operational quantities}




The order between the operational quantities derived from the norm and the injection modulus $j$ is showed in the following diagram, where  $``\rightarrow"$ means $``\leq"$:

\begin{center}
\begin{tikzpicture}[x=1cm,y=1cm]

\draw [thick,->] (0,0) -- (1.5,0);
\draw [thick,->] (2.5,0) -- (4.0,0);
\draw [thick,->] (5.0,0) -- (6.5,0);
\draw [thick,->] (7.5,0) -- (9.0,0);

\draw [] (0,0) node[left] {\normalsize $j$};
\draw [] (1.5,0) node[right] {\normalsize $\s$};
\draw [] (4.0,0) node[right] {\normalsize $\s_d$};
\draw [] (6.5,0) node[right] {\normalsize $\tau_d$};
\draw [] (9.0,0) node[right] {\normalsize $\tau$};

\draw [thick,->] (1.8,0.5) -- (1.8,2);
\draw [thick,->] (4.3,0.5) -- (4.3,2);
\draw [thick,->] (6.8,0.5) -- (6.8,2);
\draw [thick,->] (9.3,0.5) -- (9.3,2);

\draw [thick,->] (2.5,2.5) -- (4.0,2.5);
\draw [thick,->] (5.0,2.5) -- (6.5,2.5);
\draw [thick,->] (7.5,2.5) -- (9.0,2.5);
\draw [thick,->] (10,2.5) -- (11.5,2.5);

\draw [] (1.5,2.5) node[right] {\normalsize $\G$};
\draw [] (4,2.5) node[right] {\normalsize $\G_d$};
\draw [] (6.5,2.5) node[right] {\normalsize $\D_d$};
\draw [] (9,2.5) node[right] {\normalsize $\D$};
\draw [] (11.5,2.5) node[right] {\normalsize $\; \| \cdot \|$};
\end{tikzpicture}
\end{center}

The vertical arrows in the above diagram connect quantities that characterize the same classes of operators: upper semi-Fredholm, $\dns$, $\dss$ and strictly singular. We observe that none of these pairs are equivalent quantities.  

Indeed, the quantities $\s$ and $\G$ are not equivalent because $\ell_2$ is arbitrarily distortable. Hence, by \cite[Theorem 3.4 and Corollary 3.5]{GM:95}, there exist spaces $Y_n\simeq \ell_2$ and operators $T_n \in \Lc (\ell_2,Y_n)$ ($n \in \N$) such that $n \cdot \s(T_n) \leq \G(T_n)$. Since $\ell_2$ is an atomic Banach lattice, $\s_d (T_n) = \s(T_n)$ and $\G_d(T_n) = \G(T_n)$;  hence $\s_d$ and $\G_d$ are not equivalent. 



Similarly, by \cite[Proposition 1]{M:10}, the operators $T_n \in \Lc (\ell_2,Y_n)$ in the previous paragraph satisfy $n\cdot \tau(T_n) \leq \D(T_n)$,  showing that $\tau$ and $\D$ are not equivalent, and also that $\tau_d$ and $\D_d$ are not equivalent.

\subsection{Open Questions}
We finish the paper stating some open questions. 

\begin{Question}
Is $\s_d\leq D\cdot \beta$ for some constant $D>0$? 
\end{Question} 

If $E$ is an order continuous Banach lattice then $E$ is an ideal in $E^{**}$ \cite[Theorem 1.b.16]{LT-II:79}, hence the quotient $E^{**}/E$ is a Banach lattice \cite[Section  1.a]{LT-II:79}. 
Moreover, every operator $T\in\Lc(E,Y)$ induces a \emph{residuum operator}  $T^{co}\in\Lc(E^{**}/E, Y^{**}/Y)$ defined by $T^{co}(x^{**}+E)= T^{**}x^{**}+Y$.

\begin{Question}
Suppose that $E$ is  order continuous and $T\in\dns(E,Y)$. 
Is $T^{co}\in\dns$? 
\end{Question}

It was proved in \cite{GM:97} that the answer is positive in the case $E=L_1$. We refer to \cite{GSK:95} for information on the residuum operator $T^{co}$. 
\medskip

In \cite[Theorem 3.16]{GMM:20} it is shown that for $1\leq p<\infty$, $\dss(L_p,Y)$ is the perturbation class of $\dns(L_p,Y)$ in the sense that when $\dns(L_p,Y)\neq\emptyset$, $K\in\Lc(L_p,Y)$ is $\dss$ if and only if $T+K\in\dns$ for each $T\in\dns(L_p,Y)$. 

\begin{Question}
Suppose that $E$ is an order continuous Banach lattice and $\dns(E,Y)\neq\emptyset$. 

Is $\dss(E,Y)$ the perturbation class of $\dns(E,Y)$?
\end{Question}
\medskip

\noindent\textbf{Acknowledgements.}
We thank the referees for a careful reading of the manuscript and some suggestions that improved the paper.

\end{document}